\documentclass[12pt,a4paper,reqno]{amsart}
\usepackage[english]{babel}
\usepackage{amssymb,latexsym,amsfonts,amsthm,upref,amsmath}
\usepackage[margin=1in]{geometry}
\usepackage[foot]{amsaddr}
\usepackage[dvipsnames,x11names]{xcolor}
\definecolor{myblue}{HTML}{003399}
\usepackage{hyperref}
\hypersetup{colorlinks,citecolor=myblue,filecolor=black,linkcolor=myblue,urlcolor=myblue}
 
\usepackage{enumerate} 
\usepackage{tikz}
\usepackage{float}
\usepackage{cleveref}
\usepackage{mathtools}
\usepackage[noadjust]{cite}
\usepackage{filecontents}
\usepackage [autostyle, english = american]{csquotes}

\newtheorem*{corintro*}{Corollary}
\newtheorem*{thm*}{Theorem}
\newtheorem*{lem*}{Lemma}
\newtheoremstyle{prim}{}{}{\normalfont}{}{\bfseries}{.}{ }{}
\newtheoremstyle{stil}{}{}{\slshape}{}{\bfseries}{.}{ }{}
\theoremstyle{stil}
\newtheorem{thm}{Theorem}[section]
\newtheoremstyle{defi}{}{}{}{}{\bfseries}{.}{ }{}
\theoremstyle{defi}

\theoremstyle{defi}
\newtheorem{rem}[thm]{Remark}
\theoremstyle{stil}
\newtheorem{pro}[thm]{Proposition}
\theoremstyle{stil}
\newtheorem{lem}[thm]{Lemma}
\theoremstyle{stil}
\newtheorem{kor}[thm]{Corollary}
\theoremstyle{prim}

\newenvironment{prf}{\noindent \textit{Proof.}}{\null\hfill$\qed$\hskip
2mm\vskip 2mm}

\newcommand{\Y}{ {\rm Y}^+(\mathfrak{sl}_2)}

\newcommand{\R}{ {\overline{R}}}

\newcommand{\vac}{\mathop{\mathrm{\boldsymbol{1}}}}

\newcommand{\gll}{\mathfrak{gl}}

\newcommand{\sll}{\mathfrak{sl}}
\newcommand{\wsll}{\widehat{\mathfrak{sl}}}
\newcommand{\xa}{x_{\alpha}}
\newcommand{\xap}{x_{\alpha}^+}
\newcommand{\xaa}{x_{-\alpha}}

\newcommand{\vlc}{v_{N(c\Lambda_0)}}
\newcommand{\vlkn}{v_{N(k\Lambda_0)}}
\newcommand{\vlk}{v_{L(k\Lambda_0)}}

\newcommand{\CC}{\mathbb{C}}

\newcommand{\ZZ}{\mathbb{Z}}

\newcommand{\Sc}{\mathcal{S}}

\newcommand{\Tc}{\mathcal{T}}
\newcommand{\Bc}{\mathcal{B}}

\newcommand{\Vcc}{\mathcal{V}_c(\mathfrak{sl}_2)}

\newcommand{\Vck}{\mathcal{V}_k(\mathfrak{sl}_2)}

\newcommand{\Lck}{\mathcal{L}_k(\mathfrak{sl}_2)}

\newcommand{\wht}{\widehat}
\newcommand{\wvr}{\overline}

\newcommand{\ot}{\otimes}
\newcommand{\ts}{\hspace{1pt}}
\newcommand{\qdet}{ {\rm qdet}\hspace{1pt}}
\newcommand{\tr}{ {\rm tr}}

\newcommand{\ndo}{\mathop{\mathrm{End}}}
\newcommand{\om}{\mathop{\mathrm{Hom}}}

\newcommand{\cdotrl}{\mathop{\hspace{-2pt}\underset{\text{RL}}{\cdot}\hspace{-2pt}}}
\newcommand{\cdotlr}{\mathop{\hspace{-2pt}\underset{\text{LR}}{\cdot}\hspace{-2pt}}}

\newcommand{\fand}{\quad\text{and}\quad}
\newcommand{\Fand}{\qquad\text{and}\qquad}

\newcommand{\non}{\nonumber}
\newcommand{\beq}{\begin{equation}}
\newcommand{\eeq}{\end{equation}}
\newcommand{\ben}{\begin{equation*}}
\newcommand{\een}{\end{equation*}}

\makeatletter
\def\smalloverbrace#1{\mathop{\vbox{\m@th\ialign{##\crcr\noalign{\kern3\p@}%
  \tiny\downbracefill\crcr\noalign{\kern3\p@\nointerlineskip}%
  $\hfil\displaystyle{#1}\hfil$\crcr}}}\limits}
\makeatother

\makeatletter
\def\smallunderbrace#1{\mathop{\vtop{\m@th\ialign{##\crcr
   $\hfil\displaystyle{#1}\hfil$\crcr
   \noalign{\kern3\p@\nointerlineskip}%
   \tiny\upbracefill\crcr\noalign{\kern3\p@}}}}\limits}
\makeatother

\makeatletter
\def\author@andify{%
  \nxandlist {\unskip ,\penalty-1 \space\ignorespaces}%
    {\unskip {} \@@and~}%
    {\unskip \penalty-2 \space \@@and~}%
}
\makeatother

\begin{document}

\title{Principal subspaces for the quantum affine vertex algebra in type $A_1^{(1)}$}

\author{Marijana Butorac}
\address[M. Butorac]{Department of Mathematics, University of Rijeka, Radmile Matej\v{c}i\'{c} 2, 51000 Rijeka, Croatia}
\email{mbutorac@math.uniri.hr}

\author{Slaven Ko\v{z}i\'{c}} 
\address[S. Ko\v{z}i\'{c}]{ Department of Mathematics, Faculty of Science, University of Zagreb, 10000 Zagreb, Croatia}
\email{kslaven@math.hr}



\begin{abstract}
By using the ideas of Feigin and Stoyanovsky and Calinescu, Lepowsky and Milas  
we introduce and study the principal subspaces associated with the Etingof--Kazhdan quantum affine vertex algebra of integer level $k\geqslant 1$ and type $A_1^{(1)}$. We show that the principal subspaces possess the quantum vertex algebra structure, which turns   to the usual vertex algebra structure of the principal subspaces of generalized Verma and standard modules at the classical limit. Moreover,  we find their topological quasi-particle bases which  correspond  to the sum sides of certain Rogers--Ramanujan-type identities. 
\end{abstract}

\maketitle

\allowdisplaybreaks

\section{Introduction}
\numberwithin{equation}{section}

  Feigin and Stoyanovsky \cite{FS} introduced the notion of principal subspace  which can be associated with   generalized Verma modules and  integrable highest weight modules of   affine Kac--Moody Lie algebras. These rather remarkable objects 
provide a  connection between representation theory of affine  Lie algebras and combinatorial identities via their quasi-particle bases; see, e.g., \cite{FS,G,Bu1,BK}. Many other interesting aspects of principal subspaces were also extensively studied; see, e.g., \cite{CapLM,Kan,P,S} and references therein.

In this paper  we consider the Etingof--Kazhdan quantum affine vertex algebra $\Vcc$  associated with the rational $R$-matrix, as defined in \cite{EK5}. Motivated by the aforementioned results on the principal subspaces for affine Lie algebras we
  investigate their quantum counterparts using the structure of  $\Vcc$. By  generalizing the definition  of Feigin and Stoyanovsky \cite{FS} we introduce the principal subspace $W_{\Vcc}$ for $\Vcc$ and  show  that it is a quantum vertex subalgebra of $\Vcc$ for all $c\in\CC$. 
Although it  turns to the  commutative vertex algebra at the classical limit, $W_{\Vcc}$ by itself is not commutative and, furthermore, its vertex operator  map $Y(z)$ possesses   poles at $z=0$ of infinite order, i.e. $Y(v,z)w$ can have infinitely many negative powers of the variable $z$ for $v,w\in W_{\Vcc}$.  

Motivated by the presentations of principal subspaces found by Calinescu, Lepowsky and Milas \cite{CLM1,CLM2}, we introduce and study a family of   ideals $\mathcal{I}_{ \Vck }^t$, $t\in\CC$, of the quantum vertex algebra $W_{\Vck}$, where  $k>0$ is an integer. The form of the ideal generators for $t=0$ comes from the integrability relation,
	$\xa^+ (z)^{k+1}=0$ on the level $k$ standard $\wht{\sll}_2$-module, found by Lepowsky and Primc \cite{LP}. On the other hand, the generators for $\mathcal{I}_{ \Vck }^t$ with $t\neq 0$  are inspired by the $h$-adic integrability relation of vertex operators in Iohara's  bosonic realization  of level $1$ modules for the double Yangian for $\sll_2$ \cite{I}, $x^+(z)x^+(z+h)=0$.
We show that the classical limit of the  quotient quantum vertex algebra
  $W_{\Lck}=W_{\Vck}/\mathcal{I}_{ \Vck }^0$    coincides with the vertex algebra over the principal subspace of the level $k$ integrable highest weight $\wht{\sll}_2$-module $L(k\Lambda_0)$.
	In particular, using the underlying quantum vertex algebra structure, we define  quantum quasi-particles,  certain operators on
$W_{\Vcc}$ and $W_{\Lck}$, and construct    topological   bases for  $W_{\Vcc}$ and $W_{\Lck}$ consisting of certain quasi-particle monomials   applied on the vacuum vector. 
	As with  the corresponding quasi-particle bases in the affine Lie algebra case   \cite{Bu1,FS,G}, they   give an  interpretation of the sum sides of certain combinatorial identities. In the end, we further discuss quantum quasi-particle relations.

\section{Principal subspaces of certain \texorpdfstring{$\wsll_2$}{sl2hat}-modules}
Consider the complex simple Lie algebra $\sll_2
=\CC x_{-\alpha}\oplus \CC h\oplus\CC x_{\alpha}
$, where
$$
\xa =\begin{pmatrix}0&1 \\ 0&0 \end{pmatrix},\quad
h=\begin{pmatrix}1&0\\0&-1\end{pmatrix},\quad
\xaa =\begin{pmatrix}0&0\\1&0\end{pmatrix}
$$
and the Lie bracket is given by
$
[h,\xa]=2\xa$, $
[h,\xaa]=-2\xaa$ and $
[\xa,\xaa]=h$.
It is equipped with the standard symmetric invariant bilinear form,  $\left<a,b\right>=\tr (ab)$ for $a,b\in\sll_2$; see    \cite{Hum} for more details.
Let $\wsll_2=\sll_2\ot \CC[t^{\pm 1}]\oplus\CC K$ be the corresponding affine Kac--Moody Lie algebra. Its generators $K$ and $x(r)=x\ot t^r$ with $x\in\sll_2$ and $r\in\ZZ$ are subject to the relations
$$
[a(m),b(n)]=[a,b](m+n)+m\left<a,b\right>\delta_{m+n\ts 0}\ts K\fand [K,a(m)]=0.
$$ 
Fix $c\in\CC$. The generalized Verma module $N(c\Lambda_0)$ for $\wsll_2$ is defined as the induced module
$$
N(c\Lambda_0) = U(\wsll_2)\ot_{U(\sll_2\ot\CC[t] \oplus \CC K)} \CC_c ,
$$
where $U(\mathfrak{g})$ is the universal enveloping algebra of the Lie algebra $\mathfrak{g}$ and the structure of $U(\sll_2\ot\CC[t] \oplus \CC K)$-module on the one-dimensional space $\CC_c=\CC$ is defined so that $\sll_2\ot\CC[t]$ acts trivially and $K$ acts as the scalar $c$. We refer to $\vlc=1\ot 1\in N(c\Lambda_0)$ as the highest weight vector.
By the Poincar\'{e}–Birkhoff–Witt theorem the generalized Verma module $N(c\Lambda_0)$ is isomorphic to  $U(\sll_2\ot t^{-1}\CC[t^{-1}] )$ as a vector space.  
For any integer $k>0$  let $L(k\Lambda_0)$  be the level $k$ standard module, i.e.  the integrable highest weight $\wsll_2$-module which equals the unique simple quotient of  $N(k\Lambda_0)$. Let $\vlk$ be its highest weight vector, i.e. the image of the   $\vlkn $ in $L(k\Lambda_0)$; see \cite{Kac} for more details.
 It is well-known that the generalized Verma module $N(c\Lambda_0)$ and the   standard module $ L(k\Lambda_0)$ are naturally equipped with the structure of vertex algebra; see \cite{Bor,FLM,FZ}.

The principal subspaces were introduced in 
 \cite{FS}. For a generalized Verma module $V= N(c\Lambda_0)$   and   standard module $V= L(k\Lambda_0)$  define its principal subspace $W_V$ as
\beq\label{fsdef}
W_V = U(\CC \xa\ot\CC[t^{\pm 1}])\cdot v_V \subset V.
\eeq
 They can be described in terms of their presentations or in terms of  their quasi-particle bases. 
We now  recall  Calinescu--Lepowsky--Milas'  presentations      \cite{CLM1,CLM2}. 
Set
\begin{align}
&I_{N(k\Lambda_0)}=\sum_{p\geqslant k+1} U^{-}\ts R_{{N(k\Lambda_0)}}(p)
 \,\,\text{ for }\,\,  R_{{N(k\Lambda_0)}}(p) =\sum_{\substack{r_1,\ldots ,r_{k+1}\leqslant -1\\ r_1+\ldots +r_{k+1}=-p}} \xa (r_1) \ldots\xa (r_{k+1}),\label{clmp}
\end{align}
where $U^{-}=U(\CC \xa\ot t^{-1}\CC[t^{-1}])$ and $k>0$ is an integer.
Note that
 $R_{{N(k\Lambda_0)}}(p) \cdot \vlkn$   equals    the coefficient of $u^{p-k-1}$ in the   series $x_\alpha(u)^{k+1} \vlkn =x^+_\alpha(u)^{k+1} \vlkn$, where 
$$x^+_\alpha(u)=\sum_{r=1}^\infty x_\alpha(-r) u^{r-1}\fand x_\alpha(u)=\sum_{r\in\ZZ} x_\alpha(-r) u^{r-1}.$$
For any positive integer $k$ let
$ \mathcal{I}_{N(k\Lambda_0)}=I_{N(k\Lambda_0)} \cdot \vlkn$.
Recall that  the principal subspace $W_{N(c\Lambda_0)}$ with $c\in\CC$ is a vertex subalgebra of $N(c\Lambda_0)$ with the vertex operator map
\begin{align}
&Y(\xap (u_1)\ldots \xap (u_n)\vlc, z )\ts \xap (v_1)\ldots \xap (v_m)\vlc\non\\
&\qquad=\xap (z+ u_1)\ldots \xap (z+ u_n)\ts\xap (v_1)\ldots \xap (v_m)\vlc .\label{cymap}
\end{align}

\begin{thm}[\cite{CLM1,CLM2}] \label{clmthm}
$ \mathcal{I}_{N(k\Lambda_0)} $ is the   ideal of the vertex algebra $W_{N(k\Lambda_0)}$ generated by the element $\xa(-1)^{k+1}\cdot \vlkn$. 
Moreover, the  quotient vertex algebra $W_{N(k\Lambda_0)} / \mathcal{I}_{N(k\Lambda_0)} $ coincides with  the  vertex algebra $W_{L(k\Lambda_0)}$.
\end{thm}

We now recall the construction of  quasi-particle bases for principal subspaces; see \cite{Bu1,FS,G}. For any  integer $m>0$ set
\beq\label{qp1}
x_{m\alpha}^+ (u)=\sum_{r\geqslant m} x_{m\alpha}  (-r)u^{r-m}=x_{\alpha}^+(u)^m.
\eeq
The coefficients  $x_{m\alpha}(-r)$ are called quasi-particles of charge $m$. Consider the   monomials
\beq\label{mons}
x_{m_r\alpha}(n_r)\ldots x_{m_1\alpha}(n_1)\quad\text{with}\quad r\geqslant 0,\,m_1\geqslant \ldots \geqslant m_r \geqslant 1 ,\,  n_1,\ldots,n_r\leqslant -1.
\eeq
Denote by $B_{N(c\Lambda_0)} $ the set of all quasi-particle monomials \eqref{mons} which satisfy
\beq\label{cond}
n_{s+1}\leqslant n_s -2m_s \quad\text{if}\quad m_{s+1}=m_s\Fand 
n_{s}\leqslant -m_s-2(s-1) m_s
\eeq
for all $s=1,\ldots ,r-1$. The set $$\Bc_{N(c\Lambda_0)}= \big\{  b\cdot \vlc\,:\, b\in B_{N(c\Lambda_0)}\big\}\subset W_{N(c\Lambda_0)}$$ forms a basis for $W_{N(c\Lambda_0)}$; see \cite{Bu1, G}.
As for the standard modules, by  \cite{FS,G}  the set 
\begin{gather*}
\Bc_{L(k\Lambda_0)}= \left\{b\cdot \vlk\,:\, b\in B_{L(k\Lambda_0)}\right\}\subset W_{L(k\Lambda_0)},\qquad\text{where}\\ 
B_{L(k\Lambda_0)}=
\left\{
b\in  B_{N(k\Lambda_0)} \,:\, b\text{ consists of quasi-particles of charges less than or equal to }k
\right\},
\end{gather*}
forms a basis for the principal subspace $W_{L(k\Lambda_0)}$ for all integers $k>0$. Note that in the   definition of the set $\Bc_{L(k\Lambda_0)}$   the quasi-particle monomials are regarded  as operators on $L(k\Lambda_0)$.   
The bases $\Bc_{L(k\Lambda_0)}$ provide an interpretation of certain Rogers--Ramanujan-type identities, while   $\Bc_{N(c\Lambda_0)}$ can be combined with   the Poincar\'{e}–Birkhoff–Witt bases to derive new combinatorial identities; see \cite[Sect. 7]{BK} for more details and references.

\section{Quantum affine vertex algebra in type  \texorpdfstring{$A_1^{(1)}$}{A1(1)}}\label{sec02}

Let $h$ be a formal parameter.
In this paper we employ the usual expansion convention where the negative powers of the expressions of the form $z_1 +\ldots +z_n$, with each $z_i$ denoting a variable or a parameter, are expanded in the negative powers of the variable on the left. For example, we have
$$
(z_1+z_2)^{-r}=\sum_{l\geqslant 0} \binom{-r}{l} z_1^{-r-l} z_2^l
\neq 
\sum_{l\geqslant 0} \binom{-r}{l} z_2^{-r-l} z_1^l
= (z_2+z_1)^{-r}
\quad\text{for } r>0.
$$ 

  Define the Yang $R$-matrix in the variable $u$  by
\beq\label{yang}
R(u)=1-hPu^{-1}\in\ndo\CC^2\ot \ndo\CC^2 [h/u],
\eeq
where $1$ is the identity and $P$ the permutation operator in $\CC^2 \ot \CC^2$.
 There exists a unique series $g(u)\in 1+\textstyle\frac{h}{u}\CC[h/u]$ satisfying
$$
g(u+2h)=g(u)(1-h^2 u^{-2}).
$$
The normalized $R$-matrix $\overline{R}(u)= g(u)R(u)$ possesses the    unitarity property, 
\beq\label{uni}
\wvr{R} (u)\ts \wvr{R} (-u)=1 
\eeq
 and
the   crossing symmetry  properties, 
\beq\label{csym}
\wvr{R} (-u) \cdotrl\wvr{R} (u+2h )=1
\fand
\wvr{R} (-u) \cdotlr\wvr{R} (u+2h )=1,
\eeq
where the subscript RL (LR) indicates that the first tensor factor of $\wvr{R} (-u)=\wvr{R}_{12}(-u) $ is applied from the right (left) while its second tensor factor   is applied from the left (right); see, e.g., \cite[Sect. 2.2]{JKMY} for more details.

 The dual Yangian $\Y$ is the associative algebra over the commutative ring $\CC[[h]]$ with generators $t_{ij}^{(-r)}$, where $i,j=1,2$ and $r=1,2,\ldots ,$ and the defining relations
\begin{gather}
R(u-v)\ts T_1^+(u)\ts T_2^+(v)=
T_2^+(v)\ts T_1^+(u)\ts R(u-v),\label{rtt} \\
\qdet T^{+}(u) =1.\label{qdet}
\end{gather}
Regarding the  relations \eqref{rtt}, the element $T^+(u)\in \ndo\CC^2 \ot \Y[[u]]$ is defined by
\begin{align*}
T^{+}(u)=\sum_{i,j=1,2} e_{ij}\ot t_{ij}^{+}(u) \quad\text{for}\quad 
 t_{ij}^+ (u)=\delta_{ij}-h\sum_{r=1}^{\infty} t_{ij}^{(-r)}u^{r-1},
\end{align*}
where $e_{ij}$ are the matrix units.
Throughout the paper, we  indicate a copy of the matrix in the tensor product algebra
$ (\ndo\CC^2)^{\ot m}\ot\Y$ by subscripts,
so that, e.g., we have
$$
T_r^+(u)=\sum_{i,j=1,2} 1^{\otimes (r-1)} \ot e_{ij} \ot 1^{\ot (m-r)} \ot t_{ij}^+(u).
$$
In particular, we have $m=2$ and $r=1,2$ in the defining relation \eqref{rtt}. As for the relation \eqref{qdet}, for more information on the quantum determinant see \cite[Sect. 3]{EK3}.

Let us recall the Etingof--Kazhdan construction \cite[Sect. 2]{EK5} of the  quantum vertex algebra structure  over  the $h$-adically completed dual Yangian; see \cite[Sect. 1.4]{EK5} for the precise definition of the notion of quantum vertex algebra (quantum VOA).
From now on the tensor products are understood as $h$-adically completed.
 For families of  variables $u=(u_1,\ldots ,u_n)$ and $v=(v_1,\ldots ,v_m)$  and  a single variable $z$ define
$$
T_{[n]}^+(u|z)=\prod_{i=1,\dots,n}^{\longrightarrow}T_i^+(z+u_i) \fand
\R_{nm}^{12}(u|v|z)= \prod_{i=1,\dots,n}^{\longrightarrow} 
\prod_{j=n+1,\ldots,n+m}^{\longleftarrow} \R_{ij}(z+u_i -v_{j-n}) ,
$$
where the arrows indicate the order of the factors. Relation \eqref{rtt} implies the identity
\beq\label{rttgen2}
\R_{nm}^{12}(u|v|z-w)\ts T_{[n]}^{+13}(u|z)\ts T_{[m]}^{+23}(u|w)
=
T_{[m]}^{+23}(u|w)\ts T_{[n]}^{+13}(u|z)\ts\R_{nm}^{12}(u|v|z-w)
\eeq
for the operators on\vspace{-5pt}
\beq\label{ops}
\smalloverbrace{(\ndo\mathbb{C}^{2})^{\otimes n}}^{1} \otimes
\smalloverbrace{(\ndo\mathbb{C}^{2})^{\otimes m}}^{2}\otimes 
\smalloverbrace{\Y[[h]]}^{3},
\eeq
where the superscripts in \eqref{rttgen2} indicate the tensor factors, as given by \eqref{ops}.
Omitting the variable $z$ we get the operators,
$$
T_{[n]}^+(u)=\prod_{i=1,\dots,n}^{\longrightarrow}T_i^+(u_i) \fand
\R_{nm}^{12}(u|v)= \prod_{i=1,\dots,n}^{\longrightarrow} 
\prod_{j=n+1,\ldots,n+m}^{\longleftarrow} \R_{ij}(u_i -v_{j-n})
$$
which satisfy the identity
$$
\R_{nm}^{12}(u|v)\ts T_{[n]}^{+13}(u)\ts T_{[m]}^{+23}(u)
=
T_{[m]}^{+23}(u)\ts T_{[n]}^{+13}(u)\ts\R_{nm}^{12}(u|v)
$$
 on \eqref{ops}.
The next lemma is a quantum current reformulation of \cite[Lemma 2.1]{EK5}; cf. \cite{c11}.

\begin{lem}
Let $c\in\CC$ and set $\Vcc=\Y [[h]]$. 
For any  integer $n\geqslant 1$ and the family of variables $u=(u_1,\ldots ,u_n)$ there exists a unique operator
$$
\Tc_{[n]}(u)=
\Tc_{[n]}(u_1,\ldots ,u_n)\in(\ndo\CC^2)^{\ot n} \ot \om (\Vcc,\Vcc((u_1,\ldots ,u_n))[[h]])
$$
such that for all $m\geqslant 0$ and the family of variables $v=(v_1,\ldots ,v_m)$ we have
\beq\label{tmap}
\Tc_{[n]}^{13}(u)\ts T_{[m]}^{+23}(v)\vac=
\R_{nm}^{12}(u+h(c+2)|v)\cdotrl\left(T_{[n+m]}^{+}(u,v)\vac\ts \R_{nm}^{12}(u|v)^{-1}\right)
\eeq
on
\eqref{ops},
where $\vac$ is the unit in the algebra $\Y $ and $u\pm ha =(u_1 \pm ha,\ldots ,u_n\pm ha)$ for $a\in\CC$.
In particular, we have $\Tc_{[n]} (u)\vac=T_{[n]}^{+ }(u)\vac$.
\end{lem}

The next theorem is a particular case of  Etingof--Kazhdan's construction \cite[Thm. 2.3]{EK5}. 
\begin{thm}\label{EK:qva}
For any $c\in \CC$
there exists a unique  structure of quantum vertex algebra
on  $\Vcc$ such that the vacuum vector is
$\vac $, the vertex operator map is defined by
\beq\label{ymap}
Y\big(T_{[n]}^+ (u)\vac ,z\big)=\Tc_{[n]}(u|z) ,
\eeq
where $\Tc_{[n]}(u|z)=\Tc_{[n]}(z+u_1,\ldots ,z+u_n)$,
  and the braiding map $\mathcal{S}(z)$ is defined by  
\begin{align}
&\mathcal{S}(z)\big(  T_{[n]}^{+13}(u) \ts T_{[m]}^{+24}(v)   (\vac\ot  \vac) \big)
= \R_{nm}^{  12}(u|v|z-h  (c+2))^{-1}\cdotlr \Big( \R_{nm}^{  12}(u|v|z)\ts
\non\\
 &\qquad \times T_{[n]}^{+13}(u) \ts \R_{nm}^{  12}(u|v|z+h  c)^{-1} \ts
 T_{[m]}^{+24}(v)  \ts\R_{nm}^{  12}(u|v|z)(\vac\ot  \vac)\Big)\label{smap}
\end{align}
for operators on
$
(\ndo\mathbb{C}^{2})^{\otimes n} \otimes
(\ndo\mathbb{C}^{2})^{\otimes m}\otimes \Vcc \ot \Vcc $.\footnote{The original expression for the braiding map in \cite{EK5} differs from \eqref{smap}. However, by using  the crossing symmetry property \eqref{csym} one easily checks that both   definitions coincide.}
\end{thm}

\section{Principal subspaces for the quantum vertex algebra \texorpdfstring{$\Vcc$}{Vc(sl2)}}

\subsection{Principal subspaces  of  vacuum modules}

Let $c\in\CC$.
From now on  we write $x(-r)=t_{12}^{(-r)}\in\Y $ for   $r=1,2,\dots $ and 
\beq\label{evop}
x^+(u)=\sum_{r=1}^\infty x(-r) u^{r-1} \in  \Y[[u]].
\eeq
Defining relations \eqref{rtt} for the dual Yangian imply
$$
[x^+(u),x^+(v)]=0
.$$
 Let $W_{\Vcc}$ be the $h$-adic completion of the unital subalgebra of $\Y$ generated by all  $x(-r)$ with $r=1,2,\dots.$
We adopt the  terminology coming from the representation theory of  affine Lie algebras    and refer to $W_{\Vcc}$ as the principal subspace.
However, we should mention that, unlike the principal subspaces for affine Lie algebras,  $W_{\Vcc}$ is not a complex vector space but rather a topologically free $\CC[[h]]$-module.
By the Poincar\'{e}–Birkhoff–Witt theorem for the dual Yangian, see \cite[Sect. 3.4]{EK3}, the set 
$$
\Bc_{\Vcc}=\big\{x(r_m)\ldots x(r_1)\vac\,:\, r_m\leqslant \ldots \leqslant r_1\leqslant -1,\, m\geqslant 0\big\}
$$
forms a topological basis for $W_{\Vcc}$. By this we mean that  $\Bc_{\Vcc}$ is linearly independent over   $\CC[[h]]$ and that for every  $n\geqslant 1$ and $v\in W_{\Vcc}$ there exists an element $w$ in the $\CC[[h]]$-span of $\Bc_{\Vcc}$ such that
$v=w\mod h^n$.
The classical limit of $W_{\Vcc}$, i.e. the complex algebra $W_{\Vcc}/h W_{\Vcc}$ is the universal enveloping algebra $U(\CC \xa\ot t^{- 1}\CC[t^{- 1}])$. Moreover, the images of the generators $x(-r)$ with $r\geqslant 1$ in the quotient  $W_{\Vcc}/h W_{\Vcc}$ coincide with the elements  $\xa (-r)=\xa \ot t^{-r} \in U(\CC \xa\ot t^{- 1}\CC[t^{- 1}])$.

\begin{rem}\label{rem1}
As discussed above,   the universal enveloping algebra $U(\CC \xa\ot t^{- 1}\CC[t^{- 1}])$ is the classical limit  of the  algebra $W_{\Vcc}$. More generally, by taking the classical limit $h=0$ of the  dual Yangian for $ \sll_2$ one obtains the algebra $U(\sll_2\ot t^{-1}\CC[t^{-1}])$. This correspondence, as well as its generalization to the $\sll_N$  case, extends to the underlying (quantum) vertex algebra structures; see \cite[Sect. 2.2]{EK5}. Another way to go  from the quantum to the classical setting  can be found in, e.g.,  \cite[Sect. 2.2]{JKMY} and \cite[Sect. 15]{Naz} in terms of the dual (and also double) Yangian for $\gll_N$ defined over the complex field $\CC$. The classical limit $U(\gll_N\ot t^{-1}\CC[t^{-1}])$ is there obtained   as the corresponding graded algebra of the dual Yangian with respect to a certain ascending filtration.
\end{rem}

Motivated by the  quasi-particles for affine Lie algebras, see \cite{FS,G,Bu1}, we  define the series
\beq\label{qp2}
x_{(m)}^{+,t} (u)=\sum_{r\geqslant m} x_{(m)}^{t}  (-r)u^{r-m}=x^+(u)\ts x^+(u+th)\ldots x^+(u+(m-1)th)
\eeq
for any integer  $m>0$ and $t\in \CC$. In particular, we have
$x_{(m)}^{+,0} (u)=x^{+} (u)^m$.
The coefficients   $x_{(m)}^{t}  (-r)$ for $m>1$ and $t\neq 0$ are no longer elements of the dual Yangian $\Y$ because of the $h$-shifted arguments. However, they do belong to its $h$-adic completion and, more specifically, to $W_{\Vcc}$. Observe that  $x_{(1)}^{+,t} (u)=x^+ (u)$ for all $t\in\CC$. Following \cite{FS,G,Bu1} we refer to the coefficients $x_{(m)}^t (-r)$ as quasi-particles of charge $m$. 
Consider the set $$ \Bc_{\Vcc}^t=\textstyle\big\{\textstyle b\cdot \vac\,:\, b\in B_{\Vcc}^t\big\} ,$$ where
$B_{\Vcc}^t$ is the set of all monomials
$$
x_{(m_r )}^t(n_r)\ldots x_{(m_1 )}^t(n_1)\quad\text{with}\quad r\geqslant 0,\,m_1\geqslant \ldots \geqslant m_r \geqslant 1 ,\,  n_1,\ldots,n_r\leqslant -1.
$$
which satisfy difference conditions \eqref{cond} for all $s=1,\ldots ,r-1$.

\begin{pro}\label{verma_baza}
For all $c,t\in\CC$ the set $\Bc_{\Vcc}^t$ forms a topological basis for $W_{\Vcc}$.
\end{pro}

\begin{prf}
By comparing \eqref{qp1} and \eqref{qp2} we see that the classical limit of the quasi-particle  $x_{(m)}^t (-r)$ coincides with the quasi-particle $x_{m\alpha}(-r)$.
Hence the linear independence of the set $\Bc_{\Vcc}^t$ is clear as its classical limit is   the basis $\Bc_{N(c\Lambda_0)}$ of the principal subspace $W_{N(c\Lambda_0)}$.
Choose any  nonzero $v\in W_{\Vcc}$ and positive integer $n$. We will prove that there exists an element $w$ in the $\CC[[h]]$-span of $\Bc_{\Vcc}^t$ such that $v=w\mod h^n$. As $W_{\Vcc}$ is $h$-adically complete, this verifies the proposition. 
Let $\bar{v}\in W_{\Vcc}/hW_{\Vcc} =   W_{N(c\Lambda_0)}$ be the image of $v$ with respect to the classical limit. 
The $\CC[[h]]$-module $W_{\Vcc}$ is topologically free so we can assume that $\bar{v}$ is nonzero. 
Indeed, otherwise we consider the element $v'=h^{-m}v$ instead of $v$, where the integer $m>0$ is chosen so that  $v'$  belongs to $ W_{\Vcc}$ and possesses nonzero classical limit.
There exist basis elements $\bar{b}_{1,0},\ldots, \bar{b}_{r_0,0} \in \Bc_{N(c\Lambda_0)}$ and  $\beta_{1,0},\ldots, \beta_{r_0,0}\in \CC$ such that $\bar{v}=\sum_i \beta_{i,0}\ts \bar{b}_{i,0}$. Choose $b_{1,0},\ldots, b_{r_0,0}\in\Bc_{\Vcc}^t$ such that their classical limits   equal  $\bar{b}_{1,0},\ldots, \bar{b}_{r_0,0} $ respectively. This implies $v=\sum_i \beta_{i,0} \ts b_{i,0}\mod h$. Next, we express the classical limit $\bar{v}_1$ of $v_1=h^{-1}\left(v-\sum \beta_{i,0}\ts b_{i,0}\right)$ as a linear combination $\bar{v}_1=\sum_i \beta_{i,1} \ts \bar{b}_{i,1}$ for some $\bar{b}_{1,1},\ldots, \bar{b}_{r_1,1}\in \Bc_{N(c\Lambda_0)}$ and   $\beta_{1,1},\ldots, \beta_{r_1,1}\in \CC$. This produces the equality 
$v_1=\sum_i \beta_{i,1} \ts b_{i,1}\mod h$, where    $b_{i,1}\in \Bc_{\Vcc}^t$ are chosen so that their classical limits coincide with the corresponding elements $\bar{b}_{i,1}$. Therefore, we can express $v$ as
$$
v= \textstyle\sum_i \beta_{i,0}\ts b_{i,0}+h \sum_i \beta_{i,1} b_{i,1}\mod h^2.
$$
By continuing this procedure, now starting with $v_2=h^{-1}\left(v_1 -\sum_i \beta_{i,1} \ts b_{i,1}\right)$, after $n-2$ more steps we obtain a $\CC[h]$-linear combination of some  elements of $\Bc_{\Vcc}^t$ which coincides with $v$ modulo $h^n$, as required.
\end{prf}

We now show that $W_{\Vcc}$ inherits the quantum vertex algebra structure from $\Vcc$.
Let  $z$, $u=(u_1,\ldots ,u_n)$ and $v=(v_1,\ldots ,v_m)$ be the variables. Introduce the functions     
\beq\label{fje}
 g_{nm}(u|v|z)= \prod_{\substack{i=1,\dots,n\\j=1,\ldots,m}} 
   g(z+u_i -v_{j}) 
\fand 
p_{nm}(u|v|z)=
\prod_{\substack{i=1,\dots,n\\j=1,\ldots,m}} \left(1-\frac{h}{z+u_i-v_{j}}\right).
\eeq
Furthermore, we write
$x^+_{[n]}(u)=x^+(u_1)\ldots x^+ (u_n)$ and $z+u=(z+u_1,\ldots ,z+u_n)$.

\begin{thm}
$W_{\Vcc}$ is a quantum vertex subalgebra of $\Vcc$. Moreover, for all integers $m,n\geqslant 1$ and the variables $u=(u_1,\ldots ,u_n)$ and $v=(v_1,\ldots ,v_m)$ we have
\begin{align}
&Y(x^+_{[n]}(u)\vac, z)\ts x^+_{[m]}(v)\vac = p_{nm}(-u|-v|-z) \ts\non\\
& \qquad\times g_{nm}(-u|-v|-z) \ts g_{nm}(u+h(c+2)|v|z)\ts   x^+_{[n+m]}(z+u,v)\vac , \label{yfor}\\
& \Sc(z) \big(  x^+_{[n]}(u)\vac \ot \ts x^+_{[m]}(v)\vac\big)
 =p_{nm}(u|v|z)^2\ts g_{nm}(u|v|z)^2 \ts g_{nm}(-u|-v|-z-hc)\non\\
&\qquad\times g_{nm}(-u|-v|-z+(c+2)h)\ts \big(  x^+_{[n]}(u)\vac \ot \ts x^+_{[m]}(v)\vac\big).\label{sfor}
\end{align}
\end{thm}

\begin{prf}
As the vacuum vector $\vac$ belongs to $ W_{\Vcc}$ it is sufficient to prove that 
\beq\label{stst1}
Y(a,z)b\in W_{\Vcc}[[z^{\pm 1}]]\fand  \Sc(z) a\ot b \in W_{\Vcc}\ot W_{\Vcc} [[z^{\pm 1}]]\quad 
\eeq
for all 
$a,b\in W_{\Vcc}$.
In fact, it is sufficient to verify \eqref{stst1}  for the  basis elements $a,b\in \Bc_{\Vcc}$ only. However, this     follows directly  from formulae \eqref{yfor} and \eqref{sfor}.
As for the aforementioned formulae, by combining \eqref{tmap} and \eqref{ymap} we find
\begin{align}
&Y(T_{[n]}^{+13}(u)\vac,z)\ts T_{[m]}^{+23}(v)\vac\non\\
&\qquad =\R_{nm}^{12}(u+h(c+2)|v|z)\cdotrl\left(T_{[n+m]}^{+}(z+u,v)\vac\ts \R_{nm}^{12}(u|v|z)^{-1}\right).\label{yymap}
\end{align}
By using the unitarity property \eqref{uni} and taking the matrix entries of $e_{12}^{\ot (n+m)}$ in \eqref{yymap} and \eqref{smap} we obtain the identities \eqref{yfor} and \eqref{sfor} respectively, as required.
\end{prf}

\begin{rem}\label{rem2}
Due to the form of the functions in \eqref{fje}, by taking the classical limit $h=0$ of \eqref{yfor} we obtain the vertex operator map \eqref{cymap}, while the classical limit of the braiding map  \eqref{sfor} is the identity.\footnote{The same conclusion follows from the more general results in \cite{EK5}; also recall Remark \ref{rem1}.} Although the vertex operator map  for $N(c\Lambda_0)$ is commutative, this is no longer true for the quantum vertex algebra structure on $\Vcc$.  Its vertex operator map possesses  the $\Sc$-locality property, as defined in \cite[Sect. 1.3]{EK5}. For example, the particular case of the $\Sc$-locality for the vertex operator $Y(x(-1)\vac,z)$ takes the form
\begin{align}
&g(z_1-z_2)\ts g(-z_1+z_2-hc)\left(1-\frac{h}{z_1 -z_2}\right)
Y(x(-1)\vac,z_1)\ts Y(x(-1)\vac,z_2)\non\\
&\qquad =
g(z_2 -z_1)\ts g(-z_2+z_1-hc)\left(1-\frac{h}{z_2 -z_1}\right)
Y(x(-1)\vac,z_2)\ts Y(x(-1)\vac,z_1).\label{qcrell}
\end{align}
In addition,     by \eqref{yfor} the vertex operator map possesses poles  of infinite order at $z=0$.
It is worth noting that, due to  the identity  $g(z)g(-z+2h)=1$ (see \cite[Sect. 2.2]{JKMY}), equality \eqref{qcrell} at the critical level  $c=-2$ produces one of the   defining relations for the double Yangian for $\sll_2$ (see \cite[Cor. 3.4]{I}),
$$
(z_1 -z_2-h ) E_1(z_1)
\ts E_1 ( z_2)\non\\
 =
(z_1 -z_2+h )
E_1 ( z_2)\ts E_1(z_1)\quad\text{with}\quad E_1(z) =Y(x(-1)\vac,z).
$$
\end{rem}

 \subsection{Ideals of  principal subspaces}
Let $V$ be a quantum vertex algebra with the vertex operator map $Y(z)$ and the braiding map $\Sc(z)$. One can introduce the notion of ideal of $V$ in parallel with vertex algebra theory.
First, as in \cite[Def. 3.4]{Li}, for any   $I\subseteq V$ define
$$[I]=\left\{v\in V\,:\, h^n v\in I\text{ for some }n\geqslant 0\right\}. $$
A topologically free $\CC[[h]]$-submodule $I$ of $V$ 
is said to be an  ideal  of $V$
if $
I=[I]$ and
\begin{align}
&Y(z) (V\ot I),\, Y(z) (I\ot V)\subseteq I((z))[[h]] ,\label{assum2}\\
&\Sc(z) (V\ot I)\subseteq V\ot I((z))[[h]]\fand \Sc(z) (I\ot V)\subseteq I\ot V((z))[[h]].\label{assum22}
\end{align}
As with    vertex algebra theory, one  easily verifies   that the quotient of a quantum vertex algebra over its ideal is naturally equipped with the quantum vertex algebra structure. In particular, the assumption  that $I$ is a topologically free $\CC[[h]]$-module satisfying $
I=[I]$  ensures that  the quotient $V/I$ is   topologically free.

Let $k>0$ be an  integer and $t\in\CC$. 
Write $R_{\Vck}^{t}(p)=x^{t}_{(k+1)}(-p)$ for $p> k$ so  that $R_{\Vck}^t (p)\cdot\vac$ is the coefficient of $u^{p-k-1}$ in $x_{(k+1)}^{+,t} (u)\cdot \vac\in \Vck$.
Motivated by the Calinescu--Lepowsky--Milas presentations of principal subspaces from Theorem \ref{clmthm},  also recall  \eqref{clmp}, we define
\beq\label{com2}
\mathcal{I}_{\Vck}^t= \left[  I_{\Vck}^t \cdot \vac\right] [[h]] \quad\text{for}\quad
I_{\Vck}^t=\sum_{p\geqslant k+1} W_{\Vck} R_{\Vck}^t(p),
\eeq
where the action of $W_{\Vck}$ on $R_{\Vck}^t(p)$ is given by the algebra multiplication.

\begin{pro} \label{qpres}
$\mathcal{I}_{\Vck}^t $ is the ideal of the quantum vertex algebra $W_{\Vck}$. 
\end{pro} 

\begin{prf}
First, we note that $\mathcal{I}_{\Vck}^t$ is topologically free by  construction. Next, the constraint $\mathcal{I}_{\Vck}^t =[\mathcal{I}_{\Vck}^t ]$ follows by \cite[Prop. 3.7]{Li}. Hence it remains to verify the requirements imposed by \eqref{assum2} and \eqref{assum22}. Choose any two basis monomials $a,b\in B_{\Vck}^t$ and $p\geqslant k+1$. Let us prove that the images of
$a  R_{\Vck}^t(p) \ot b$ and $b\ot a  R_{\Vck}^t(p) $ with respect to the vertex operator map 
$Y(z)$ and braiding $\Sc(z)$ belong to 
\beq\label{belong}
\mathcal{I}_{\Vck}^t((z))[[h]]\quad\text{and to}\quad  \mathcal{I}_{\Vck}^t \ot \Vck((z))[[h]]\text{ and } \Vck\ot \mathcal{I}_{\Vck}^t((z))[[h]]
\eeq
respectively.
Apply the substitution 
$$
(u_{n-k},  \ldots ,u_n)=(w,   w+th,\ldots ,  w+kth)
\quad\text{or}\quad
(v_{n-k},  \ldots ,v_n)=(w,   w+th,\ldots ,  w+kth),
$$
where $w$ is a single variable and $n >k$ or $m>k$ arbitrary integers, to formulae \eqref{yfor} and \eqref{sfor}. This turns the expressions $x^+_{[n+m]}(z+u,v)\vac$ and $x^+_{[n]}(u)\vac \ot \ts x^+_{[m]}(v)\vac$, which appear on
the right hand sides of \eqref{yfor} and \eqref{sfor}, into  power series   with coefficients in $\mathcal{I}_{\Vck}^t $  and in $\mathcal{I}_{\Vck}^t  \ot \Vck$ or $ \Vck\ot \mathcal{I}_{\Vck}^t $ respectively. Hence the same happens with the  entire right hand sides of \eqref{yfor} and \eqref{sfor} as well. 
Therefore, by extracting the  coefficients for suitably chosen $n$ and $m$ we conclude that  the images of
$a  R_{\Vck}^t(p) \ot b$ and $b\ot a  R_{\Vck}^t(p) $ with respect to  
$Y(z)$ and $\Sc(z)$ belong to \eqref{belong} respectively, as required. Finally, as  $Y(z)$ and $\Sc(z)$ are $\CC[[h]]$-module maps, this conclusion extends from the elements of the form $a  R_{\Vck}^t(p) \ot b$ and $b\ot a  R_{\Vck}^t(p) $ to all elements of the $h$-adically completed tensor products $\mathcal{I}_{\Vck}^t \ot \Vck$ and $\Vck\ot \mathcal{I}_{\Vck}^t$, so the proposition follows.
\end{prf}

Denote by $\bar{B}^{t}_{\Vck}$   the set of all monomials in $B_{\Vck}^t$ which contain at least one quasi-particle of charge greater than or equal to $k+1$. We write
$\bar{\Bc}^t_{\Vck}=\textstyle\big\{\textstyle b\cdot \vac\,:\, b\in \bar{B}^{t}_{\Vck}\big\}$. By Proposition \ref{verma_baza} the  set $\bar{\Bc}^t_{\Vck}$ is linearly independent.
Consider the identities
$$
x_{(k+l+1)}^{+,t}(z) = x_{(l)}^{+,t}(z+(k+1)th)\ts x^{+,t}_{(k+1)}(z)\quad\text{for }l> 0,\, t\in\CC.
$$
They imply that for any integer $n\geqslant 1$ each element of $\bar{\Bc}^t_{\Vck}$ can be expressed as a $\CC[h]$-linear combination of the elements of   $I_{\Vck}^t \cdot \vac $ modulo $h^n$. Therefore, as the ideal  $\mathcal{I}_{\Vck}^t $ is $h$-adically complete,   it contains the   set $\bar{\Bc}^t_{\Vck}$. 
From now on, we consider the $t=0$ case. As for the $t\neq 0$ case, see Section \ref{tneqo} and, in particular, Remark \ref{tnn} below.

\begin{pro}\label{ideal}
The set $\bar{\Bc}^0_{\Vck}$ forms a topological basis of  $\mathcal{I}_{\Vck}^0 $.
\end{pro}

\begin{prf}
By the discussion above the given set is linearly independent. The proof that  $\bar{\Bc}^0_{\Vck}$   spans an $h$-adically dense $\CC[[h]]$-submodule of $\mathcal{I}_{\Vck}^0 $  goes by repeating the argument of Jerkovi\'{c} and Primc \cite[Sect. 4.4]{JP}.
It relies on the fact that  quasi-particles of charges $p$ and $q$ with $p\leqslant q$ satisfy $2p$ independent relations of the form
\beq\label{relss}
\left(\frac{d^l}{dz^l} x^{+,0}_{(p)} (z)  \right)x^{+,0}_{(q)} (z) 
=A_l(z) ,\quad l=0,\ldots , 2p-1.
\eeq
Here $A_0(z),\ldots ,A_{2p-1}(z)$ are certain formal power series with coefficients in the set of quasi-particle polynomials such that each coefficient contains at least one quasi-particle of charge greater than or equal to $q+1$. By arguing as in the affine Lie algebra setting, one checks that for any integer $n\geqslant 1$   relations \eqref{relss} can be used to express
any element of $\mathcal{I}_{\Vck}^0 $ as a $\CC[h]$-linear combination of the elements of $\bar{\Bc}^0_{\Vck}$  modulo $h^n$.
\end{prf}

Due to Proposition \ref{qpres}  we can introduce the quotient quantum vertex algebra
$$W_{\Lck}= W_{\Vck} / \mathcal{I}_{\Vck}^0 .$$ 
We now regard the monomials in $B_{\Vck}^0$  as operators on $W_{\Lck}$. Define
$$
\Bc_{\Lck}= \left\{b\cdot \vac\,:\, b\in B_{\Vck}^0\setminus \bar{B}^0_{\Vck}\right\}\subset W_{\Lck}. 
$$
Clearly, the set $B_{\Vck}^0\setminus \bar{B}^0_{\Vck}$ consists of all quasi-particle monomials in $B_{\Vck}^0$  whose quasi-particle charges are less than or equal to $k$.

\begin{thm}
The classical limit of the quantum vertex algebra $W_{\Lck}$ coincides with the vertex algebra $W_{L(k\Lambda_0)}$. Moreover, the set $\Bc_{\Lck}$ 
forms a topological basis of $W_{\Lck}$.
\end{thm}

\begin{prf}
The fact that $\Bc_{\Lck}$ spans an $h$-adically dense $\CC[[h]]$-submodule of $W_{\Lck}$ goes by a  usual argument which relies on the quasi-particle relations
\eqref{relss} and closely follows Jerkovi\'{c} and Primc \cite[Sect. 4.4]{JP};  also recall the proof of  Proposition \ref{ideal}.
  Moreover, the set $\Bc_{\Lck}$ is linearly independent   as it is complementary to the topological basis $\bar{\Bc}^0_{\Vck}$ of the ideal $\mathcal{I}^0_{\Vck} $ established in Proposition \ref{ideal}. Finally, as the quantum vertex algebra structure on $W_{\Lck}$ comes from \eqref{yfor} and \eqref{sfor}, we conclude by Remark \ref{rem2} that its classical limit produces the vertex algebra $W_{L(k\Lambda_0)}$, as required.
\end{prf}

\begin{rem}
The quantum analogue
 of the  defining identity for quasi-particles, 
\beq\label{clc}
x_{m\alpha}(z)=x_{\alpha}(z)^m=Y(\xa(-1)^m \vlk,z),
\eeq   
see  \cite[Sect. 3]{G}, takes the form
\beq\label{lbl}
x_{\left<m\right>}^t(z)\coloneqq Y(x^{m,t}(-1)\vac,z),\,\text{ where }\, 
x^{m,t}(-1)=x(-1)x^+(th)\ldots x^+((m-1)th),
\eeq
  $x^+(ath)$ denote the evaluations of the series \eqref{evop} at $u=ath$ with $t\in\CC$ and $x (-1)$ is the constant term of  \eqref{evop}. 
	In particular, by setting $t=0$  in \eqref{lbl} we get
	\beq\label{pwr}
	x_{\left<m\right>}^0(z)=Y(x^{m,0}(-1)\vac,z)=Y(x(-1)^m\vac,z).
	\eeq
	However, in contrast with \eqref{clc}, the expression in \eqref{pwr} is not the $m$-th power of the series $x^0_{\left<1\right>}(z)$. In  fact, the $m$-th power of $x^0_{\left<1\right>}(z)$ is not even well defined.
	Also note that $x_{\left<m\right>}^t (z)\vac = x_{(m)}^{+,t} (z)\vac$ for all $t\in\CC$ and, furthermore, that the   classical limits of $x_{\left<m\right>}^t (z)$ and $x_{(m)}^{+,t} (z) $   coincide. On the other hand, the actions of the operators $x_{\left<m\right>}^{t}(z)$ and $x^{+,t}_{(m)}(z)$ on the quotient $ W_{\Vck} / \mathcal{I}_{\Vck}^t $ differ when $k\geqslant m$ and  are both trivial for $k<m$. 
	The requirement $x^{+,t}_{(k+1)}(z)=0$ on   level $k$, which defines the ideal $\mathcal{I}_{\Vck}^t$ in \eqref{com2}, may be regarded as a quantum version of the   integrability relation,
	$\xa^+ (z)^{k+1}=0$ on the level $k$ standard $\wht{\sll}_2$-module of Lepowsky and Primc \cite{LP}. 
	Also, in view of Remark
	\ref{rem2}, it is worth noting that the operator $E_1 (z)$, which  corresponds to $x_{\left<1\right>}(z)=Y(x(-1)\vac, z)$ and comes from   Iohara's bosonic realization \cite[Thm. 4.7]{I} of level $1$ modules for   the double Yangian for $\sll_2$, satisfies the analogous constraint, $E_1(z)E_1(z+h)=0$ on level $1$. 
\end{rem}

As with the Calinescu--Lepowsky--Milas presentation \cite{CLM1,CLM2}, the ideal   $\mathcal{I}_{\Vck}^0 $ turns out to be principal; recall Theorem \ref{clmthm}.

\begin{kor}
For any integer $k\geqslant 1$  the ideal  $\mathcal{I}_{\Vck}^0 $  of the quantum vertex algebra $ W_{\Lck}$ is generated by the element $x(-1)^{k+1} \vac$.
\end{kor}

\begin{prf}
Let $J$ be an ideal of $ W_{\Lck}$ which contains $x(-1)^{k+1} \vac$.
First, we observe  that 
$x(-1)^{k+1}\vac $ is the constant term of $x_{(k+1)}^{+,0} (u)\vac$, i.e. we have
$x(-1)^{k+1}\vac = R_{\Vck}^0 (k+1)\vac$, so that $x(-1)^{k+1} \vac $ belongs to  
$\mathcal{I}_{\Vck}^0 $ as well.
Next, by 
$$Y(x(-1)^{k+1}\vac,z)\vac= x^{+,0}_{(k+1)}(z)\vac =\sum_{p\geqslant k+1} R_{\Vck}^0(p)\vac z^{p-k-1}$$
we conclude that  all elements  $R^0_{\Vck}(p)\vac$ with $p\geqslant k+1$ belong to $J$. Consider the components  $x_{\left<1\right>}^0(r)$   of the vertex operator
$$
Y(x(-1)\vac, z)=x_{\left<1\right>}^0(z)=\sum_{r\in\ZZ} x_{\left<1\right>}^0(r) z^{-r-1}.
$$
Identity    \eqref{yfor} implies that the action of all $x(s)$, $s\leqslant -1$,
can be expressed using $x_{\left<1\right>}^0(r)$,  $r\in\ZZ$. Therefore, we have $\mathcal{I}_{\Vck}^0 \subseteq J$  so the corollary  follows.
\end{prf}

 \subsection{Quasi-particle relations at \texorpdfstring{$t\neq 0$}{t!=0}}\label{tneqo}

In contrast with  \cite{FS,G,Bu1}, where   the construction of  quasi-particle bases of principal subspaces relies on finding a suitable family of quasi-particle relations, in the proof of Proposition \ref{verma_baza} we bypassed this step using the classical limit. However, the quantum quasi-particles, as defined by   \eqref{qp2}, still satisfy certain relations which can be employed instead of classical limit to directly verify that the   set    $\Bc_{\Vcc}^t$ spans  $h$-adically dense $\CC[[h]]$-submodule  of the corresponding principal subspace.
For  $t=0$ these are given by   \eqref{relss}.
Let us derive such relations for $t\neq 0$.
First, for integers $q\geqslant p>0$ we  have  the following $2p$ identities:
\begin{align}
&
x^{+,t}_{(p)} (z+c_k h) \ts x^{+,t}_{(q)} (z+pth)  = x^{+,t}_{(k-1)} (z+pth)\ts  x^{+,t}_{(p+q-k+1)} (z+c_k h),\quad k=1,\ldots ,p,\label{rhs1}\\
&
 x^{+,t}_{(p)} (z+c_{p+k} h) \ts x^{+,t}_{(q)} (z+pth)  = x^{+,t}_{(p-k)} (z+c_{p+k} h)\ts  x^{+,t}_{( q+k)} (z+pth),\quad k=1,\ldots ,p,\label{rhs2}
\end{align}
where  we set $x^{+,t}_{(0)} (z)=1$ and 
$$(c_1,\ldots ,c_p ;c_{p+1},\ldots ,c_{2p})=(0,t,\ldots,t(p-1); t(q+1),t(q+2),\ldots ,t(q+p)).$$ 

\begin{rem}
Relations \eqref{rhs1} and \eqref{rhs2} present a rational $R$-matrix counterpart of the quasi-particle relations for the Ding--Feigin operators associated with the quantum affine algebra in type $A$; cf. \cite[Sect. 3.1]{c01}.
\end{rem}

 For   $l=1,\ldots ,2p$ let $(\alpha_{1,l},\ldots , \alpha_{l,l})$ be the solution of the system of $l$ linear equations
\beq\label{system}
c_1^k \alpha_{1,l} + c_2^k \alpha_{2,l} +\ldots + c_l^k \alpha_{l,l} = \delta_{k\ts l-1} (l-1)!, \quad k=0,\ldots ,l-1,
\eeq
where $c_1^k =\delta_{k\ts 0}$. Clearly,  each of these $2p$ systems possesses  a unique solution  as its coefficients form a Vandermonde matrix. Indeed, we have   $c_i\neq c_j$ for $i\neq j$ because of $q\geqslant p$.
Denote the expressions on the right hand sides of  equalities \eqref{rhs1} and \eqref{rhs2} obtained for $k=1,\ldots ,p$ by $R_k^t$ and $R_{p+k}^t$ respectively. We have $2p$  relations,
\beq\label{rhs3}
(th)^{1-l}\sum_{k=1}^l \alpha_{k,l}\ts x_{(p)}^{+,t} (z+c_k h) \ts x_{(q)}^{+,t} (z+p h)
 =(th)^{1-l} \sum_{k=1}^l \alpha_{k,l} R_k^t\quad\text{for} \quad l=1,\ldots ,2p.
\eeq
By using \eqref{system} and the formal Taylor theorem, 
$a(z+z_0)=e^{z_0 \frac{d}{dz}} a(z)
 $ 
one   checks that the sum on the left hand side of \eqref{rhs3} possesses a zero of order greater than or equal to $l-1$ at $h=0$  so that the   expressions in \eqref{rhs3} are well defined.
Relations \eqref{rhs3} are  linear combinations of equalities \eqref{rhs1} and \eqref{rhs2} and    each coefficient on the right hand side of \eqref{rhs3}   contains one quasi-particle of charge strictly greater than $q$.
 Furthermore, by \eqref{system} we have
\beq\label{rhs4}
(th)^{1-l}\sum_{k=1}^l \alpha_{k,l}\ts x_{(p)}^{+,t} (z+c_k h) \ts x_{(q)}^{+,t} (z+p h)=
\left( \frac{d^{l-1}}{dz^{l-1}}\ts x_{(p)}^{+,t} (z)\right) x_{(q)}^{+,t} (z )\mod h 
\eeq
for $l=1,\ldots ,2p$.
As the classical limit of the right hand side of \eqref{rhs4} is equal to 
$$\left( \frac{d^{l-1}}{dz^{l-1}}\ts x_{p\alpha}^+ (z)\right) x_{q\alpha}^+ (z ),$$
it follows by   \cite[Sect. 4.4]{JP}, see also  \cite[Lemma 2.2.1]{Bu1}, that the quasi-particle relations  obtained by equating the right hand sides of \eqref{rhs3} and \eqref{rhs4},  
\beq\label{krj}
\left( \frac{d^{l-1}}{dz^{l-1}}\ts x_{(p)}^{+,t} (z)\right) x_{(q)}^{+,t} (z )
=(th)^{1-l} \sum_{k=1}^l \alpha_{k,l} R_k^t
\mod h
\qquad\text{for }   l=1,\ldots ,2p
\eeq
  are  independent. More specifically, they can be used in parallel with the corresponding quasi-particle relations  in the affine Lie algebra setting  \cite{FS,G,Bu1}   to  prove directly that  the set   $\Bc_{\Vcc}^t$ spans  $h$-adically dense $\CC[[h]]$-submodule  of the corresponding principal subspace. However,  as  relations \eqref{krj} are given modulo $h$, the classical argument has to be applied $n$ times   to express a given element of $W_{\Vcc}$ as a $\CC[h]$-linear combination of elements of $\Bc_{\Vcc}^t$ modulo $h^n$. Consequently, the argument relies on the fact that the principal subspace $V=  W_{\Vcc}$ is topologically free and  possesses the following property:
	\beq\label{krj2}
	a=b\mod h\quad\text{for  }a,b\in V\qquad\text{implies}\qquad  h^{-1}(a-b) \in V .
	\eeq

	\begin{rem}\label{tnn}
	It is not clear whether the constraint in \eqref{krj2} holds for $V= I_{\Vck}^t \cdot \vac$. Hence the aforementioned classical argument does not necessarily lead to the construction of topological basis of $\mathcal{I}_{\Vck}^t$ with $t\neq 0$.
A major difference between the ideal $\mathcal{I}_{\Vck}^0$ and the ideals $\mathcal{I}_{\Vck}^t$ with $t\neq 0$ appears to be the fact that $\mathcal{I}_{\Vck}^0$ is of the form $ \mathcal{I}_{N(k\Lambda_0)}[[h]] $. 
On the other hand, consider the following simple example.
Let $t\neq 0$ and $k=1$. By regarding $x_{(2)}^{+,t}(z)\vac\in\mathcal{I}_{\mathcal{V}_1(\sll_2)}^t $ as a power series with respect to the parameter $h$ and  then taking the classical limit of each coefficient we get
$$
 x_{(2)}^{+,t}(z)   \vac
\equiv x_{2\alpha} (z)v_{N(\Lambda_0)} +h\frac{t}{2}\frac{d}{dz}x_{2\alpha} (z)v_{N(\Lambda_0)} + h^2\frac{t^2 }{2}x_{\alpha} (z) \frac{d^2}{dz^2}x_{\alpha} (z)v_{N(\Lambda_0)}+\ldots  .
$$
The coefficients with respect to the variable $z$ of the first two summands on the right hand side belong to $\mathcal{I}_{N(\Lambda_0)}[h] $, which is no longer true for the third summand. 
\end{rem}

\section*{Acknowledgement}
The authors are grateful to   Mirko Primc for helpful discussions.
This work has been supported in part by Croatian Science Foundation under the project UIP-2019-04-8488. The first  author is partially supported by the QuantiXLie Centre of Excellence, a project cofinanced by the Croatian Government and European Union through the European Regional Development Fund - the Competitiveness and Cohesion Operational Programme (Grant KK.01.1.1.01.0004).

\linespread{1.0}

\end{document}